\documentclass[border=10pt]{amsart}
\usepackage{amsfonts,amsmath, amscd, amssymb, yhmath,mathtools}
\usepackage{pgf,tikz,verbatim}
\usepackage[all,import,rotate]{xy}
\usepackage{latexsym,mathrsfs,url}
\xyoption{all}

\newtheorem{theorem}{Theorem}[section]

\newtheorem{lem}[theorem]{Lemma}

\newtheorem*{B-principle}{B-principle}
\newtheorem*{H-principle}{H-principle}

\theoremstyle{remark}

\theoremstyle{definition}

\begin{document}

\title{Enhanced Kauffman bracket}
\author{Zhiqing Yang}
\thanks{The author is supported by a grant (No. 11271058) of NSFC.}

\address{Department of Mathematics, Dalian University of Technology, China}
\email{yangzhq@dlut.edu.cn}

\date{\today}

\begin{abstract}
S. Nelson, M. Orrison, V. Rivera {\cite{S}} modified Kauffman's construction of bracket. Their invariant $\Phi^{\beta}_X$ takes value in a finite ring $Z_2[t]/(1+t+t^3)$. In this paper, the author generalizes this invariant. The new invariant takes value in a polynomial ring. Furthermore, for a tricolorable link diagram, the author gives a bracket invariant which gives lower bound on number of crossings with  different (same) colors.
\end{abstract}

\keywords{knot invariant \and knot polynomial \and Kauffman bracket }
\subjclass[2000]{Primary 57M27; Secondary 57M25}

\maketitle

\setcounter{tocdepth}{2}

\section{Introduction}\label{sec:1}

Polynomial invariants of links have a long history. In 1928, J.W. Alexander {\cite{Al}} discovered the famous Alexander polynomial. It has many connections with other topological invariants. In 1984 Vaughan Jones {\cite{J}}  discovered the Jones polynomial. Later, Louis Kauffman in 1987 {\cite{K}} introduced Kauffman bracket. It satisfies $<
\tikz{\draw[thin][scale=.06] (0,0) -- (1.5,1.5);
\draw[thin][scale=.06]  (2.5,2.5) -- (4,4);
\draw[thin][scale=.06]  (0,4) -- (4,0);} > =
A
<\tikz{\draw[thin][scale=.06] (12,12) .. controls (14,10) and (14,10) .. (12,8);
\draw[thin][scale=.06] (16,12) .. controls (14,10) and (14,10) .. (16,8);;} >
+A^{-1}
<\tikz{\draw[thin][scale=.06] (6,12) .. controls (8,10) and (8,10) .. (10,12);
\draw[thin][scale=.06] (6,8) .. controls (8,10) and (8,10) .. (10,8);} >$.

The Kauffman bracket can be calculated in two ways. First, it can be calculated inductively. In this point of view, \tikz{\draw[thin][scale=.06] (0,0) -- (1.5,1.5);
\draw[thin][scale=.06]  (2.5,2.5) -- (4,4);
\draw[thin][scale=.06]  (0,4) -- (4,0);} ,
\tikz{\draw[thin][scale=.06] (12,12) .. controls (14,10) and (14,10) .. (12,8);
\draw[thin][scale=.06] (16,12) .. controls (14,10) and (14,10) .. (16,8);;} , \tikz{\draw[thin][scale=.06] (6,12) .. controls (8,10) and (8,10) .. (10,12);
\draw[thin][scale=.06] (6,8) .. controls (8,10) and (8,10) .. (10,8);}  are regarded as three link diagrams which are identical except in a small disk. Then the calculation of $<
\tikz{\draw[thin][scale=.06] (0,0) -- (1.5,1.5);
\draw[thin][scale=.06]  (2.5,2.5) -- (4,4);
\draw[thin][scale=.06]  (0,4) -- (4,0);} >$ reduces to the calculation of $<\tikz{\draw[thin][scale=.06] (12,12) .. controls (14,10) and (14,10) .. (12,8);
\draw[thin][scale=.06] (16,12) .. controls (14,10) and (14,10) .. (16,8);;} >$ and $
<\tikz{\draw[thin][scale=.06] (6,12) .. controls (8,10) and (8,10) .. (10,12);
\draw[thin][scale=.06] (6,8) .. controls (8,10) and (8,10) .. (10,8);} >$. Second, the Kauffman bracket satisfies the following state sum formula.

$$<L>= \sum_S A^{a(S)}A^{-b(S)}(-A^2-A^{-2})^{|S|-1}$$
One can simultaneously apply $<
\tikz{\draw[thin][scale=.06] (0,0) -- (1.5,1.5);
\draw[thin][scale=.06]  (2.5,2.5) -- (4,4);
\draw[thin][scale=.06]  (0,4) -- (4,0);} > =
A
<\tikz{\draw[thin][scale=.06] (12,12) .. controls (14,10) and (14,10) .. (12,8);
\draw[thin][scale=.06] (16,12) .. controls (14,10) and (14,10) .. (16,8);;} >
+A^{-1}
<\tikz{\draw[thin][scale=.06] (6,12) .. controls (8,10) and (8,10) .. (10,12);
\draw[thin][scale=.06] (6,8) .. controls (8,10) and (8,10) .. (10,8);} >$ to all crossing points and get $2^n$ states. The righthand side is a summation over all states.

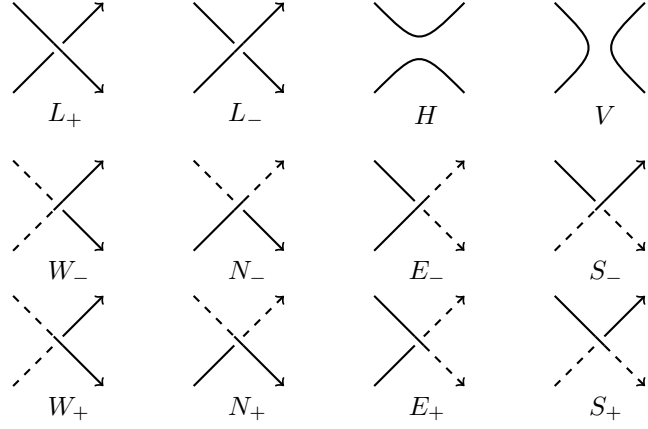
\begin{figure}[ht]
\beginpgfgraphicnamed{graphic-of-flat-world}
\begin{tikzpicture}[scale=.3]
\draw[thick][dashed] (0,0) -- (1.8,1.8);
\draw[thick] [->] (2.2,2.2) -- (4,4);
\draw[thick][dashed]  (0,4) -- (1.8,2.2);
\draw[thick] (1.8,2.2) -- (2.2,1.8);
\draw[thick] [->] (2.2,1.8) -- (4,0);
\draw (2,-1) node[text width=8] {$W_+$};

\draw[shift={(8,0)}][thick](0,0) -- (1.8,1.8);
\draw[shift={(8,0)}][thick] [dashed] [->] (2.2,2.2) -- (4,4);
\draw[shift={(8,0)}][thick][dashed]  (0,4) -- (1.8,2.2);
\draw[shift={(8,0)}][thick] (1.8,2.2) -- (2.2,1.8);
\draw[shift={(8,0)}][thick] [->] (2.2,1.8) -- (4,0);
\draw[shift={(8,0)}] (2,-1) node[text width=8] {$N_+$};

\draw[shift={(16,0)}][thick](0,0) -- (1.8,1.8);
\draw[shift={(16,0)}][thick] [dashed] [->] (2.2,2.2) -- (4,4);
\draw[shift={(16,0)}][thick]  (0,4) -- (1.8,2.2);
\draw[shift={(16,0)}][thick] (1.8,2.2) -- (2.2,1.8);
\draw[shift={(16,0)}][thick][dashed] [->] (2.2,1.8) -- (4,0);
\draw[shift={(16,0)}] (2,-1) node[text width=8] {$E_+$};

\draw[shift={(24,0)}][thick][dashed](0,0) -- (1.8,1.8);
\draw[shift={(24,0)}][thick] [->] (2.2,2.2) -- (4,4);
\draw[shift={(24,0)}][thick]  (0,4) -- (1.8,2.2);
\draw[shift={(24,0)}][thick] (1.8,2.2) -- (2.2,1.8);
\draw[shift={(24,0)}][thick][dashed]  [->] (2.2,1.8) -- (4,0);
\draw[shift={(24,0)}] (2,-1) node[text width=8] {$S_+$};

\draw[shift={(0,6)}][thick][dashed] (0,0) -- (1.8,1.8);
\draw[shift={(0,6)}][thick] [->] (2.2,2.2) -- (4,4);
\draw[shift={(0,6)}][thick][dashed]  (0,4) -- (1.8,2.2);
\draw[shift={(0,6)}][thick] (1.8,1.8) -- (2.2,2.2);
\draw[shift={(0,6)}][thick] [->] (2.2,1.8) -- (4,0);
\draw[shift={(0,6)}] (2,-1) node[text width=8] {$W_-$};

\draw[shift={(8,6)}][thick](0,0) -- (1.8,1.8);
\draw[shift={(8,6)}][thick] [dashed] [->] (2.2,2.2) -- (4,4);
\draw[shift={(8,6)}][thick][dashed]  (0,4) -- (1.8,2.2);
\draw[shift={(8,6)}][thick]  (1.8,1.8) -- (2.2,2.2);
\draw[shift={(8,6)}][thick] [->] (2.2,1.8) -- (4,0);
\draw[shift={(8,6)}] (2,-1) node[text width=8] {$N_-$};

\draw[shift={(16,6)}][thick](0,0) -- (1.8,1.8);
\draw[shift={(16,6)}][thick] [dashed] [->] (2.2,2.2) -- (4,4);
\draw[shift={(16,6)}][thick]  (0,4) -- (1.8,2.2);
\draw[shift={(16,6)}][thick]  (1.8,1.8) -- (2.2,2.2);
\draw[shift={(16,6)}][thick][dashed] [->] (2.2,1.8) -- (4,0);
\draw[shift={(16,6)}] (2,-1) node[text width=8] {$E_-$};

\draw[shift={(24,6)}][thick][dashed](0,0) -- (1.8,1.8);
\draw[shift={(24,6)}][thick] [->] (2.2,2.2) -- (4,4);
\draw[shift={(24,6)}][thick]  (0,4) -- (1.8,2.2);
\draw[shift={(24,6)}][thick]  (1.8,1.8) -- (2.2,2.2);
\draw[shift={(24,6)}][thick][dashed]  [->] (2.2,1.8) -- (4,0);
\draw[shift={(24,6)}] (2,-1) node[text width=8] {$S_-$};

\draw[shift={(0,13)}][thick] (0,0) -- (1.8,1.8);
\draw[shift={(0,13)}][thick] [->] (2.2,2.2) -- (4,4);
\draw[shift={(0,13)}][thick]  (0,4) -- (1.8,2.2);
\draw[shift={(0,13)}][thick] (1.8,2.2) -- (2.2,1.8);
\draw[shift={(0,13)}][thick] [->] (2.2,1.8) -- (4,0);
\draw[shift={(0,13)}] (2,-1) node[text width=8] {$L_+$};

\draw[shift={(8,13)}][thick](0,0) -- (1.8,1.8);
\draw[shift={(8,13)}][thick]  [->] (2.2,2.2) -- (4,4);
\draw[shift={(8,13)}][thick](0,4) -- (1.8,2.2);
\draw[shift={(8,13)}][thick] (1.8,1.8) -- (2.2,2.2);
\draw[shift={(8,13)}][thick] [->] (2.2,1.8) -- (4,0);
\draw[shift={(8,13)}] (2,-1) node[text width=8] {$L_-$};

\draw[shift={(10,5)}][thick] (6,12) .. controls (8,10) and (8,10) .. (10,12);
\draw[shift={(10,5)}][thick] (6,8) .. controls (8,10) and (8,10) .. (10,8);
\draw[shift={(10,5)}] (8,7) node[text width=4] {$H$};
\draw[shift={(12,5)}][thick] (12,12) .. controls (14,10) and (14,10) .. (12,8);
\draw[shift={(12,5)}][thick] (16,12) .. controls (14,10) and (14,10) .. (16,8);
\draw[shift={(12,5)}] (14,7) node[text width=4] {$V$};
\end{tikzpicture}
\endpgfgraphicnamed
\caption{Local crossings.}
\label{fig1}
\end{figure}

This construction can be modified. For example, S. Nelson, M. Orrison, V. Rivera {\cite{S}} introduced the following way to enhance the bracket polynomial. A link diagram can bicolored as follows. Choose two colors, say, solid and dotted. The crossing points divide any link component into even number of segments. Pick one segment and assign one color to it. Then change to another color whenever pass one crossing point. Do this to each component, we  get a bicolor link diagram. If the link has $k$ components, then there are $2^k$ different ways of colorings. Letters $E,S,W,N$ mean the east, south,west and north directions as in usual maps,  $+$ means positive crossing, $-$ means negative crossing. Now a crossing point may have one of the four types: $N_+,N-_,S_+,S_-$. See Fig. ~\ref{fig1}. If this is a virtual link diagram, then there may have four more types: $E_+,E_-,W_+,W_-$. For example, $S_+$ means that the crossing is of positive type, and the two dotted arcs are divided by the ray from the crossing point to south.

For a bicolor link diagram, one can {\cite{S}} use the following skein relations.

\begin{flalign*}
N_+ = H+ t V, &\ \ \ \ \   N_- = H+ (1+t^2) V,  \\
S_+ = H+ t V, &\ \ \ \ \     S_- = H+ (1+t^2) V,  \\
E_+ =(1+t) H+ (t+t^2) V,  &\ \ \ \ \   E_- =t H+ V,   \\
W_+ =(1+t^2) H+  V,  &\ \ \ \ \    W_- =(t+t^2)  H+ (1+t)  V,   \\
<D\  \tikz{\draw [thick][scale=.25] (0,1) -- (0,0);
\draw [thick][scale=.25] (0,0) -- (1,0);
\draw [thick][scale=.25] (1,0) -- (1,1);}\  \tikz{\draw (0,0) circle (0.14cm);}> &=(1+t+t^2) <D>.
\end{flalign*}

S. Nelson, M. Orrison, V. Rivera {\cite{S}} modified Kauffman's construction of bracket. Their invariant $\Phi^{\beta}_X$ takes value in $Z_2[t]/(1+t+t^3)$. It can distinguish some pair of knots which the HOMFLY polynomial cannot distinguish. For example, $\Phi^\beta_X(
10_{132}) = 2t + t^2\neq 2 + t + t^2 = \Phi^\beta_X(5_1)$.

\section{Enhanced Bracket}
Given any link diagram, we choose two colors, say, solid and dotted. Then the link diagram can be colored. For a bicolor link diagram, we define the following skein relations.
\begin{equation}\label{eq1}
 \left\{
\begin{lgathered}
N_+ =a_n H+ b_n V,  \ \ \ \ \   N_- =a_n' H+ b_n' V,  \\
S_+ =a_s H+ b_s V,  \ \ \ \ \   S_- =a_s' H+ b_s' V,  \\
E_+ =a_e H+ b_e V,  \ \ \ \ \   E_- =a_e' H+ b_e' V,   \\
W_+ =a_w H+ b_w V,  \ \ \ \ \   W_- =a_w' H+ b_w' V.   \\
< D\  \tikz{\draw [thick][scale=.25] (0,1) -- (0,0);
\draw [thick][scale=.25] (0,0) -- (1,0);
\draw [thick][scale=.25] (1,0) -- (1,1);}\  \tikz{\draw (0,0) circle (0.14cm);}> =d <D>.
\end{lgathered}
\right.
\end{equation}

When we apply a skein relation to a crossing, say $N_+$, we shall get two smoothings, $H$ and $V$. We will say that $H$ is the result of type $A$ smoothing, $V$ is the result of type $B$ smoothing.

As mentioned before, the Kauffman bracket can be calculated in two ways. First, it can be calculated inductively. In this point of view, the calculation of $<
\tikz{\draw[thin][scale=.06] (0,0) -- (1.5,1.5);
\draw[thin][scale=.06]  (2.5,2.5) -- (4,4);
\draw[thin][scale=.06]  (0,4) -- (4,0);} >$ reduces to the calculation of $<\tikz{\draw[thin][scale=.06] (12,12) .. controls (14,10) and (14,10) .. (12,8);
\draw[thin][scale=.06] (16,12) .. controls (14,10) and (14,10) .. (16,8);;} >$ and $
<\tikz{\draw[thin][scale=.06] (6,12) .. controls (8,10) and (8,10) .. (10,12);
\draw[thin][scale=.06] (6,8) .. controls (8,10) and (8,10) .. (10,8);} >$. Second, the Kauffman bracket satisfies the following state sum formula.

$$<L>= \sum_S A^{a(S)}A^{-b(S)}(-A^2-A^{-2})^{|S|-1}$$
One can simultaneously apply $<
\tikz{\draw[thin][scale=.06] (0,0) -- (1.5,1.5);
\draw[thin][scale=.06]  (2.5,2.5) -- (4,4);
\draw[thin][scale=.06]  (0,4) -- (4,0);} > =
A
<\tikz{\draw[thin][scale=.06] (12,12) .. controls (14,10) and (14,10) .. (12,8);
\draw[thin][scale=.06] (16,12) .. controls (14,10) and (14,10) .. (16,8);;} >
+A^{-1}
<\tikz{\draw[thin][scale=.06] (6,12) .. controls (8,10) and (8,10) .. (10,12);
\draw[thin][scale=.06] (6,8) .. controls (8,10) and (8,10) .. (10,8);} >$ to all crossing points and get $2^n$ states. The righthand side is a summation over all states.

Strictly speaking, we should not use $``="$ in skein relations (1)-(5). $``\rightarrow "$ is more appropriate, since the bracket we defined will not satisfy $<
\tikz{\draw[thin][scale=.06] (0,0) -- (1.5,1.5);
\draw[thin][scale=.06]  (2.5,2.5) -- (4,4);
\draw[thin][scale=.06]  (0,4) -- (4,0);} > =
a_n
<\tikz{\draw[thin][scale=.06] (12,12) .. controls (14,10) and (14,10) .. (12,8);
\draw[thin][scale=.06] (16,12) .. controls (14,10) and (14,10) .. (16,8);;} >
+b_n
<\tikz{\draw[thin][scale=.06] (6,12) .. controls (8,10) and (8,10) .. (10,12);
\draw[thin][scale=.06] (6,8) .. controls (8,10) and (8,10) .. (10,8);} >$ even for the first case. This is because that the three link diagrams
\tikz{\draw[thin][scale=.06] (0,0) -- (1.5,1.5);
\draw[thin][scale=.06]  (2.5,2.5) -- (4,4);
\draw[thin][scale=.06]  (0,4) -- (4,0);},
<\tikz{\draw[thin][scale=.06] (12,12) .. controls (14,10) and (14,10) .. (12,8);
\draw[thin][scale=.06] (16,12) .. controls (14,10) and (14,10) .. (16,8);;} and  \tikz{\draw[thin][scale=.06] (6,12) .. controls (8,10) and (8,10) .. (10,12);
\draw[thin][scale=.06] (6,8) .. controls (8,10) and (8,10) .. (10,8);} do not have a canonical matched coloring. At any other crossing, since the colorings are different, one will apply different skein relations to the three links. Hence $<
\tikz{\draw[thin][scale=.06] (0,0) -- (1.5,1.5);
\draw[thin][scale=.06]  (2.5,2.5) -- (4,4);
\draw[thin][scale=.06]  (0,4) -- (4,0);} > =
a_n
<\tikz{\draw[thin][scale=.06] (12,12) .. controls (14,10) and (14,10) .. (12,8);
\draw[thin][scale=.06] (16,12) .. controls (14,10) and (14,10) .. (16,8);;} >
+b_n
<\tikz{\draw[thin][scale=.06] (6,12) .. controls (8,10) and (8,10) .. (10,12);
\draw[thin][scale=.06] (6,8) .. controls (8,10) and (8,10) .. (10,8);} >$ does not hold. This is why we cannot calculate the invariant inductively.

For our construction, one can only use the second way. Namely, first step, choose a orientation of the link and one way to color the link diagram. Second, apply the above skein relations to each crossing and get $2^n$ states. Then take summation over all states. One cannot calculate the invariant inductively since $<\tikz{\draw[thin][scale=.06] (12,12) .. controls (14,10) and (14,10) .. (12,8);
\draw[thin][scale=.06] (16,12) .. controls (14,10) and (14,10) .. (16,8);;} >$ and $
<\tikz{\draw[thin][scale=.06] (6,12) .. controls (8,10) and (8,10) .. (10,12);
\draw[thin][scale=.06] (6,8) .. controls (8,10) and (8,10) .. (10,8);} >$ do not preserve the coloring of $<
\tikz{\draw[thin][scale=.06] (0,0) -- (1.5,1.5);
\draw[thin][scale=.06]  (2.5,2.5) -- (4,4);
\draw[thin][scale=.06]  (0,4) -- (4,0);} >$.

There are some relations among those coefficients to guarantee Reidemeister move invariance. We shall discuss them one by one.

\bigskip
\subsection{Reidemeister move II}

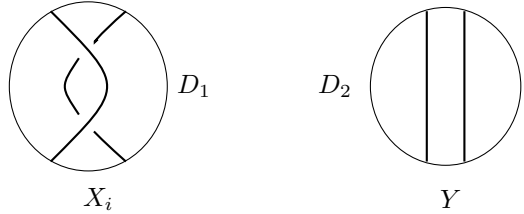
\begin{figure}[ht]
\begin{tikzpicture}[scale=.25]

\draw[thick] (0,8) .. controls (4,4) and (4,4) .. (0,0);
\draw (2,-2) node[text width=4pt] {$X_i$};
\draw[thick]  (4,8) .. controls (2.5,6.7) and (2,6)  .. (2.5,6.5);
\draw[thick]  (1.5,5.5) .. controls (0.5,4) and (0.5,4)  .. (1.5,2.5);
\draw[thick]  (4,0) .. controls (3,0.9)   .. (2.3,1.6);
 \draw (7,4) node[text width=4pt] {$D_1$};

 \draw[shift={(20,0)}][thick]  (0,0) -- (0,8);
\draw[shift={(20,0)}][thick]   (2,0) -- (2,8);
\draw (21,-2) node[text width=4pt] {$Y$};
 \draw (14.5,4) node[text width=4pt] {$D_2$};

\draw  (2,4) ellipse (4.2 and 4.5);
\draw  (21,4) ellipse (4 and 4.2);

\end{tikzpicture}
\caption{Reidemeister move II}
\label{fig2}
\end{figure}

Given two link diagrams $L_1,L_2$. Outside of the disks $D_1,D_2$, their diagrams are the same. Inside of the disks $D_1,D_2$, they are as in Fig \ref{fig2}.

\begin{figure}[ht]
\begin{tikzpicture}[scale=.25]
\draw[thick] [<-] (0,8) .. controls (4,4) and (4,4) .. (0,0); \draw (2,-1) node[text width=4pt] {$X_1$};
\draw[thick] [<-] (4,8) .. controls (2.5,6.7) and (2,6)  .. (2.5,6.5);
\draw[thick]  (1.5,5.5) .. controls (0.5,4) and (0.5,4)  .. (1.5,2.5);
\draw[thick]  (4,0) .. controls (3,0.9)   .. (2.3,1.6);
 \draw (-.1,-.1) node[text width=18pt] {$a$};
 \draw (5.5,-.1) node[text width=18pt] {$b$};

\draw[thick] [<-] (4.001,-2.8878) .. controls (0.001,-6.8878) and (0.001,-6.8878) .. (4.001,-10.8878); \draw (2.001,-11.8878) node[text width=4pt] {$X_1'$};
\draw[thick] [<-] (0.001,-2.8878) .. controls (1.501,-4.1878) and (2.001,-4.8878)  .. (1.501,-4.3878);
\draw[thick]  (2.501,-5.3878) .. controls (3.501,-6.8878) and (3.501,-6.8878)  .. (2.501,-8.3878);
\draw[thick]  (0.001,-10.8878) .. controls (1.001,-9.9878)   .. (1.701,-9.2878);
 \draw (-.1,-11.1) node[text width=18pt] {$a$};
 \draw (5.5,-11.1) node[text width=18pt] {$b$};

\draw[thick] [shift={(0,-11)}] [<-] (9.1102,19.11) .. controls (13.1102,15.11) and (13.1102,15.11) .. (9.1102,11.11);
\draw  [shift={(0,-11)}] (11.1102,10.11) node[text width=4pt] {$X_2$};
\draw[thick] [shift={(0,-11)}]  (13.1102,19.11) .. controls (11.6102,17.81) and (11.1102,17.11)  .. (11.6102,17.61);
\draw[thick] [shift={(0,-11)}] (10.6102,16.61) .. controls (9.6102,15.11) and (9.6102,15.11)  .. (10.6102,13.61);
\draw[thick] [shift={(0,-11)}] [<-] (13.1102,11.11) .. controls (12.1102,12.01)   .. (11.4102,12.71);
  \draw (9.3,-.1) node[text width=18pt] {$a$};
 \draw (14.7,-.1) node[text width=18pt] {$b$};

\draw[thick] [shift={(0,-11)}] [->] (13.1112,8) .. controls (9.1112,4) and (9.1112,4) .. (13.1112,0);
\draw  [shift={(0,-11)}] (11.1112,-1) node[text width=4pt] {$X_2'$};
\draw[thick][shift={(0,-11)}] [<-] (9.1112,8) .. controls (10.6112,6.7) and (11.1112,6)  .. (10.6112,6.5);
\draw[thick] [shift={(0,-11)}] (11.6112,5.5) .. controls (12.6112,4) and (12.6112,4)  .. (11.6112,2.5);
\draw[thick] [shift={(0,-11)}] (9.1112,0) .. controls (10.1112,0.9)   .. (10.8112,1.6);
  \draw (9.2,-11.1) node[text width=18pt] {$a$};
 \draw (14.7,-11.1) node[text width=18pt] {$b$};

\draw[thick] [shift={(0,-22)}] [->] (17.9982,29.9978) .. controls (21.9982,25.9978) and (21.9982,25.9978) .. (17.9982,21.9978);
\draw  [shift={(0,-22)}] (19.9982,20.9978) node[text width=4pt] {$X_3$};
\draw[thick] [shift={(0,-22)}] [<-] (21.9982,29.9978) .. controls (20.4982,28.6978) and (19.9982,27.9978)  .. (20.4982,28.4978);
\draw[thick] [shift={(0,-22)}] (19.4982,27.4978) .. controls (18.4982,25.9978) and (18.4982,25.9978)  .. (19.4982,24.4978);
\draw[thick] [shift={(0,-22)}]  (21.9982,21.9978) .. controls (20.9982,22.8978)   .. (20.2982,23.5978);
  \draw (18.1,-.1) node[text width=18pt] {$a$};
 \draw (23.7,-.1) node[text width=18pt] {$b$};

\draw[thick] [shift={(0,-22)}] [<-] (22.2214,18.8878) .. controls (18.2214,14.8878) and (18.2214,14.8878) .. (22.2214,10.8878);
\draw  [shift={(0,-22)}] (20.2214,9.8878) node[text width=4pt] {$X_3'$};
\draw[thick][shift={(0,-22)}]  (18.2214,18.8878) .. controls (19.7214,17.5878) and (20.2214,16.8878)  .. (19.7214,17.3878);
\draw[thick] [shift={(0,-22)}] (20.7214,16.3878) .. controls (21.7214,14.8878) and (21.7214,14.8878)  .. (20.7214,13.3878);
\draw[thick] [shift={(0,-22)}] [<-] (18.2214,10.8878) .. controls (19.2214,11.7878)   .. (19.9214,12.4878);
  \draw (18.2,-11.1) node[text width=18pt] {$a$};
 \draw (23.8,-11.1) node[text width=18pt] {$b$};
\end{tikzpicture}
\caption{Oriented Reidemeister move II}
\label{fig3}
\end{figure}
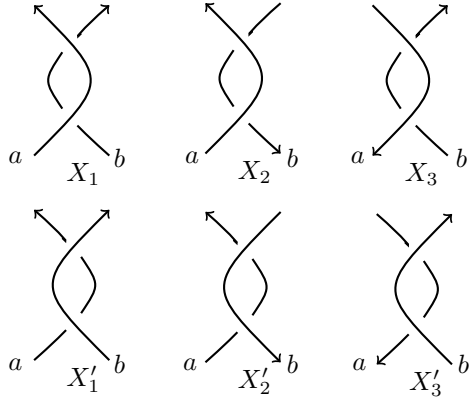

When we consider orientations, their are $6$ cases for the diagram inside $D_1$. See Fig \ref{fig3}.

To apply skein relations (1), we also need to color the link diagram. For each of the $X_i$ or $X_i'$, each of the two segments $a,b$ can be colored in 2 ways. So there are 4 cases. Let's first consider the case that they are both dotted. Then for link diagrams $L_1,L_2$, we use skein relations to all crossings outside the disks $D_1,D_2$. Let $\Lambda$ be the set of all possible smoothings
outside $D_i$, let $\Lambda_i'$ be the set of all possible smoothings
inside $D_i$. Then we have the following result.

$$<L_i>= \sum_{S\in \Lambda } f(S) \sum_{S'\in \Lambda_i' } g(S,S')$$

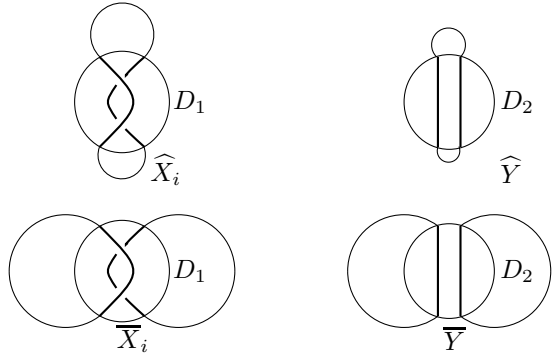
\begin{figure}[ht]
\begin{tikzpicture}[scale=.15]

\draw[thick] (0,8) .. controls (4,4) and (4,4) .. (0,0);
\draw (2,-2) node[text width=4pt] {$\overline{X}_i$};
\draw[thick]  (4,8) .. controls (2.5,6.7) and (2,6)  .. (2.5,6.5);
\draw[thick]  (1.5,5.5) .. controls (0.5,4) and (0.5,4)  .. (1.5,2.5);
\draw[thick]  (4,0) .. controls (3,0.9)   .. (2.3,1.6);
 \draw (7,4) node[text width=4pt] {$D_1$};
\draw (-8,4) arc (180:50:5);
\draw (-8,4) arc (180:310:5);
\draw (12,4) arc (0:130:5);
\draw (12,4) arc (0:-130:5);

 \draw[shift={(30,0)}][thick]  (0,0) -- (0,8);
\draw[shift={(30,0)}][thick]   (2,0) -- (2,8);
\draw (31,-2) node[text width=4pt] {$\overline{Y}$};
 \draw (36,4) node[text width=4pt] {$D_2$};

\draw  (2,4) ellipse (4.2 and 4.5);
\draw  (31,4) ellipse (4 and 4.2);

\draw (22,4) arc (180:53:5);
\draw (22,4) arc (180:307:5);
\draw (40,4) arc (0:127:5);
\draw (40,4) arc (0:-127:5);

\draw[thick][shift={(0,15)}] (0,8) .. controls (4,4) and (4,4) .. (0,0);
\draw[shift={(0,15)}] (5,-2) node[text width=4pt] {$\widehat{X}_i$};
\draw[thick][shift={(0,15)}]  (4,8) .. controls (2.5,6.7) and (2,6)  .. (2.5,6.5);
\draw[thick][shift={(0,15)}]  (1.5,5.5) .. controls (0.5,4) and (0.5,4)  .. (1.5,2.5);
\draw[thick][shift={(0,15)}]  (4,0) .. controls (3,0.9)   .. (2.3,1.6);
 \draw[shift={(0,15)}] (7,4) node[text width=4pt] {$D_1$};
\draw[shift={(0,15)}] (4,8) arc (-45:225:2.8);
\draw[shift={(0,15)}] (0,0) arc (160:383:2.1);

 \draw[shift={(30,15)}][thick]  (0,0) -- (0,8);
\draw[shift={(30,15)}][thick]   (2,0) -- (2,8);
\draw[shift={(0,15)}] (36,-2) node[text width=4pt] {$\widehat{Y}$};
 \draw[shift={(0,15)}] (36,4) node[text width=4pt] {$D_2$};

\draw[shift={(0,15)}]  (2,4) ellipse (4.2 and 4.5);
\draw[shift={(0,15)}]  (31,4) ellipse (4 and 4.2);

\draw[shift={(0,15)}] (32,8) arc (-45:225:1.5);
\draw[shift={(0,15)}] (30,0) arc (160:383:1);

\end{tikzpicture}
\caption{Outside smoothing patterns}
\label{fig4}
\end{figure}

Since there skein relations smooth the crossings, $\Lambda$ is a disjoint union of two sets. $\Lambda = \Lambda (1) \cup \Lambda (2)$, where $\Lambda (1)$ consists of all smoothings of pattern $\widehat{X}_i$ and $\widehat{Y}$ in figure {\ref{fig4}}, $\Lambda (2)$ consists of all smoothings of pattern $\overline{X}_i$ and $\overline{Y}$ in figure {\ref{fig4}}.

$$<L_i>= \sum_{S\in \Lambda (1) } f(S) \sum_{S'\in \Lambda_i' } g(S,S')+ \sum_{S\in \Lambda (2) } f(S) \sum_{S'\in \Lambda_i' } g(S,S')$$

So we have

$$<L_1>= \sum_{S\in \Lambda (1) } f(S) \sum_{S'\in \Lambda_1' } g(S,S')+ \sum_{S\in \Lambda (2) } f(S) \sum_{S'\in \Lambda_1' } g(S,S')=f_1 <\overline{X}_i>+f_2<\widehat{X}_i>$$

$$<L_2>= \sum_{S\in \Lambda (1) } f(S) \sum_{S'\in \Lambda_2' } g(S,S')+ \sum_{S\in \Lambda (2) } f(S) \sum_{S'\in \Lambda_2' } g(S,S')=f_1 <\overline{Y}>+f_2<\widehat{Y}>$$

Notice that the coefficients of $<\overline{X}_i>$ and $<\overline{Y}>$ are the same.  Hence to have $<L_1>=<L_2>$, it is sufficient to have $<\overline{X}_i>=<\overline{Y}>$ and $<\widehat{X}_i>=<\widehat{Y}>$.

In figure {\ref{fig4}}, if the color of arcs $a,b$ are dotted, then we have

$$<\overline{X}_1>= (a_wa_e')d^2+(a_wb_e'+b_wa_e'+b_wb_e'd)d$$
$$<\overline{Y}>=d^2$$
So we have $(a_wa_e')d+(a_wb_e'+b_wa_e'+b_wb_e'd)=d$. Similarly,
 $$<\widehat{X}_1>= (a_wa_e')d+(a_wb_e'+b_wa_e'+b_wb_e'd)d^2$$
$$<\widehat{Y}>=d$$
So we have $(a_wa_e')+(a_wb_e'+b_wa_e'+b_wb_e'd)d=1$.

In summary, for $X_1$, we have :

$(a_wa_e')d+(a_wb_e'+b_wa_e'+b_wb_e'd)=d$ and  $(a_wa_e')+(a_wb_e'+b_wa_e'+b_wb_e'd)d=1$.

For the linear system of equations $xd+y=d,x+yd=1$, one can get solutions $$\{ x=1,y=0\}, \{ d=1, x+y=1\}, \{ d=-1, x-y=1\}.$$ For simplicity, we just consider the case $\{ x=1,y=0\}$ here. Hence we have $$a_wa_e'=1, a_wb_e'+b_wa_e'+b_wb_e'd=0.$$

If we change the colors of $a,b$ to solid, then we will get $$a_ea_w'=1, a_eb_w'+b_ea_w'+b_eb_w'd=0.$$
If $a$ is dotted, $b$ is solid, then we will get $$a_sa_s'=1, a_sb_s'+b_sa_s'+b_sb_s'd=0.$$
If $b$ is dotted, $a$ is solid, then we will get $$a_na_n'=1, a_nb_n'+b_na_n'+b_nb_n'd=0.$$

Similarly, for $X_2,X_2'$, we will get the following set of equations.
\[
 \left\{
\begin{gathered}
b_wb_e'=1, b_wa_e'+a_wb_e'+a_wa_e'd=0 \\
b_eb_w'=1, b_ea_w'+a_eb_w'+a_ea_w'd=0\\
b_nb_n'=1, b_na_n'+a_nb_n'+a_na_n'd=0\\
b_sb_s'=1, b_sa_s'+a_sb_s'+a_sa_s'd=0.
\end{gathered}
\right.
\]

For $X_3,X_3'$, we will get the same set of equations as $X_2,X_2'$.

In summary, for Reidemeister move II, we get the following equations.

\begin{equation}\label{R2}
 \left\{
\begin{lgathered}
a_wa_e'=1,  a_wb_e'+b_wa_e'+b_wb_e'd=0\\
a_ea_w'=1, a_eb_w'+b_ea_w'+b_eb_w'd=0\\
a_sa_s'=1, a_sb_s'+b_sa_s'+b_sb_s'd=0\\
a_na_n'=1, a_nb_n'+b_na_n'+b_nb_n'd=0\\
b_wb_e'=1, b_wa_e'+a_wb_e'+a_wa_e'd=0 \\
b_eb_w'=1, b_ea_w'+a_eb_w'+a_ea_w'd=0\\
b_nb_n'=1, b_na_n'+a_nb_n'+a_na_n'd=0\\
b_sb_s'=1, b_sa_s'+a_sb_s'+a_sa_s'd=0.
\end{lgathered}
\right.
\end{equation}

\bigskip
\subsection{Reidemeister move III}

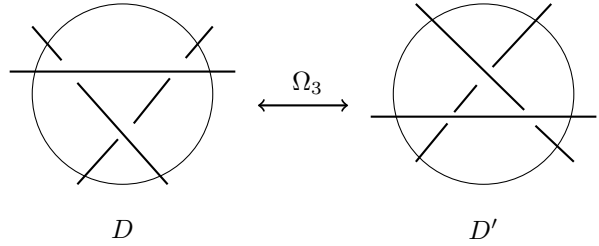
\begin{figure}[h]
\begin{tikzpicture}[scale=0.3,mine/.style={font=\itshape\large}]
\draw[thick]  (9,6) -- (19,6);
\draw[thick]   (10,8) -- (11.3,6.5);
\draw[thick]  (12,5.5) -- (16,1);

\draw[thick](18,8) -- (16.8,6.6);
\draw[thick] (16.1,5.7) -- (14.5,3.7);
\draw[thick] (13.8,3) -- (12,1);

\draw[thick] (14,-1) node {$D$};

\draw[thick][shift={(0,-.5)}]   [<->]   (20,5) -- (24,5);
\draw[thick][shift={(0,-.5)}]  (22.2,6) node {$\Omega_{3}$};

\draw[thick][shift={(0,1)}] (25,3) -- (35,3);
\draw[thick][shift={(0,1)}]  (27,8) -- (31.8,3.3);
\draw[thick][shift={(0,1)}]    (32.3,2.6) -- (34,1);
\draw[thick][shift={(0,1)}] (33,8) -- (30.4,5.2);
\draw[thick][shift={(0,1)}]  (29.7,4.4) -- (28.7,3.2);
\draw[thick][shift={(0,1)}]  (28.4,2.7) -- (27,1);
\draw[thick][shift={(0,0)}]  (30,-1) node {$D'$};

\draw  (14,5) ellipse (4 and 4);
\draw[shift={(0,1)}]   (30,4) ellipse (4 and 4);
\end{tikzpicture}

\caption{Reidemeister move three.}
\label{fig5}
\end{figure}

Given two link diagrams $L,L'$. Outside of the disks $D,D'$, their diagrams are the same. Inside of the disks $D,D'$, they are as in Fig \ref{fig5}.

\begin{figure}
\begin{tikzpicture}[scale=0.2]
\draw  (0,0) ellipse (4 and 4);
\draw (-1.5,3.6) arc (215:-35:1.9);
\draw (-4,0) arc (120:300:2.2);
\draw (4,0) arc (60:-120:2.2);
\draw [thick] (-1,-9) node[text width=0.4pt] {$Out (1)$};
\draw [dashed] (-1.5,3.6) -- (1.5,-3.6);
\draw [dashed] (1.5,3.6) -- (-1.5,-3.6);
\draw [dashed] (-4,0) -- (4,0);
\draw [dashed][shift={(14,0)}]  (-1.5,3.6) -- (1.5,-3.6);
\draw [dashed][shift={(14,0)}]  (1.5,3.6) -- (-1.5,-3.6);
\draw [dashed][shift={(14,0)}] (-4,0) -- (4,0);
\draw [dashed] [shift={(28,0)}](-1.5,3.6) -- (1.5,-3.6);
\draw [dashed] [shift={(28,0)}](1.5,3.6) -- (-1.5,-3.6);
\draw [dashed] [shift={(28,0)}](-4,0) -- (4,0);
\draw [dashed] [shift={(42,0)}] (-1.5,3.6) -- (1.5,-3.6);
\draw [dashed] [shift={(42,0)}] (1.5,3.6) -- (-1.5,-3.6);
\draw [dashed] [shift={(42,0)}] (-4,0) -- (4,0);
\draw [dashed] [shift={(56,0)}](-1.5,3.6) -- (1.5,-3.6);
\draw [dashed] [shift={(56,0)}](1.5,3.6) -- (-1.5,-3.6);
\draw [dashed] [shift={(56,0)}](-4,0) -- (4,0);

\draw [shift={(14,0)}]  (0,0) ellipse (4 and 4);
\draw [shift={(14,0)}] (-4,0) arc (240:55:2.2);
\draw [shift={(14,0)}] (4,0) arc (-60:125:2.2);
\draw [shift={(14,0)}] (-1.5,-3.6) arc (145:390:1.9);
\draw [shift={(14,0)}][thick] (-1,-9) node[text width=0.4pt] {$Out (2)$};

\draw [shift={(28,0)}]  (0,0) ellipse (4 and 4);
\draw [shift={(28,0)}] (-1.5,3.6) arc (215:-35:2);
\draw [shift={(28,0)}] (-4,0) arc (220:-42:5.3);
\draw [shift={(28,0)}] (-1.5,-3.6) arc (145:390:1.9);
\draw [shift={(28,0)}][thick] (-1,-9) node[text width=0.4pt] {$Out (3)$};

\draw [shift={(42,0)}]  (0,0) ellipse (4 and 4);
\draw [shift={(42,0)}] (-4,0) arc (240:60:2.2);
\draw [shift={(42,0)}] (1.5,3.6) arc (100:-145:4.7);
\draw [shift={(42,0)}] (4,0) arc (60:-120:2.2);
\draw [shift={(42,0)}][thick] (-1,-9) node[text width=0.4pt] {$Out (4)$};

\draw [shift={(56,0)}]  (0,0) ellipse (4 and 4);
\draw [shift={(56,0)}] (-4,0) arc (120:300:2.2);
\draw [shift={(56,0)}] (4,0) arc (-60:120:2.2);
\draw [shift={(56,0)}] (-1.5,3.6) arc (140:-100:4.5);
\draw [shift={(56,0)}][thick] (-1,-9) node[text width=0.4pt] {$Out (5)$};
\end{tikzpicture}
\caption{States outside the disks.}
\label{fig6}
\end{figure}
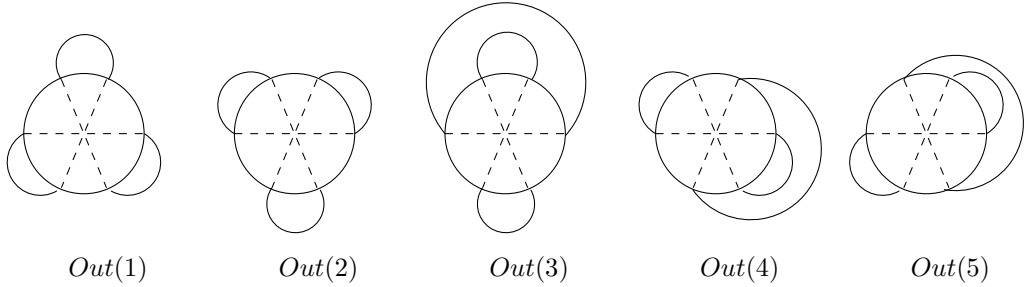

When we use skein relations to all crossings outside the disks $D,D'$, there are five patterns $Out (1)$-- $Out(5)$ as in figure {\ref{fig6}}.

Let $D_i$ be the diagram that inside a disk it is the same as $L$ in $D$, outside the disk it is the same as $Out (i)$. Let $D_i'$ be the diagram that inside a disk it is the same as $L'$ in $D'$, outside the disk it is the same as $Out (i)$.  If we smooth all outside crossings, we will have the following result.

$$<L>= \sum_{i=1}^5 f_i <D_i>,\  <L'>= \sum_{i=1}^5 f_i <D_i'>.$$

Hence $<D_i>=<D_i'>$ for all $i$ implies that $ <L>= <L'>$.

\begin{figure}[h]
\begin{tikzpicture}[scale=0.2,mine/.style={font=\itshape\large}]
\draw  (0,0) ellipse (4 and 4);
\draw [thick] (-1.5,3.6) arc (-215:35:1.9);
\draw [thick] (-4,0) arc (120:-60:2);
\draw [thick] (4,0) arc (60:240:2);
\draw [thick] (0,-5) node[text width=0.4pt] {$A$};

\draw [shift={(10,0)}]  (0,0) ellipse (4 and 4);
\draw [shift={(10,0)}][thick] (-1.5,3.6) arc (-215:35:1.9);
\draw [shift={(10,0)}][thick] (-4,0) arc (120:-60:2);
\draw [shift={(10,0)}][thick] (4,0) arc (60:240:2);
\draw [shift={(10,0)}][thick] (0,-5) node[text width=0.4pt] {$B$};
\draw  [thick] (10,-.5) ellipse (.8 and .8);

\draw [shift={(20,0)}]  (0,0) ellipse (4 and 4);
\draw [shift={(20,0)}][thick] (-4,0) arc (-120:60:2);
\draw [shift={(20,0)}][thick] (4,0) arc (300:120:2);
\draw [shift={(20,0)}][thick] (-1.5,-3.6) arc (215:-35:1.9);
\draw [shift={(20,0)}][thick] (0,-5) node[text width=0.4pt] {$C$};

\draw [shift={(30,0)}]  (0,0) ellipse (4 and 4);
\draw [shift={(30,0)}][thick] (-4,0) arc (-120:60:2);
\draw [shift={(30,0)}][thick] (4,0) arc (300:120:2);
\draw [shift={(30,0)}][thick] (-1.5,-3.6) arc (215:-35:1.9);
\draw [shift={(30,0)}][thick] (0,-5) node[text width=0.4pt] {$D$};
\draw  [thick] (30,.5) ellipse (.8 and .8);

\draw [shift={(40,0)}]  (0,0) ellipse (4 and 4);
\draw [shift={(40,0)}][thick] (-1.5,3.6) arc (-215:35:1.9);
\draw [shift={(40,0)}][thick] (-1.5,-3.6) arc (215:-35:1.9);
\draw [shift={(40,0)}][thick] (-4,0) -- (4,0);
\draw [shift={(40,0)}][thick] (0,-5) node[text width=0.4pt] {$E$};

\draw [shift={(50,0)}]  (0,0) ellipse (4 and 4);
\draw [shift={(50,0)}][thick] (4,0) arc (300:120:2);
\draw [shift={(50,0)}][thick] (-4,0) arc (120:-60:2);
\draw [shift={(50,0)}][thick]  (-1.5,3.6) -- (1.5,-3.6);
\draw [shift={(50,0)}][thick] (0,-5) node[text width=0.4pt] {$F$};

\draw [shift={(60,0)}]  (0,0) ellipse (4 and 4);
\draw [shift={(60,0)}][thick] (-4,0) arc (-120:60:2);
\draw [shift={(60,0)}][thick] (4,0) arc (60:240:2);
\draw [shift={(60,0)}][thick]  (1.5,3.6) -- (-1.5,-3.6);
\draw [shift={(60,0)}][thick] (0,-5) node[text width=0.4pt] {$G$};
\end{tikzpicture}

\caption{States inside the disks.}
\label{fig7}
\end{figure}
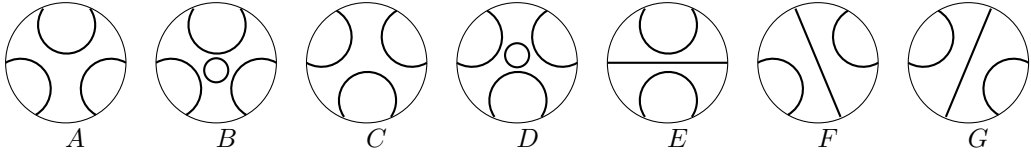

When we use skein relations to all crossings inside the disks $D,D'$, there are seven patterns $A$-- $G$ as in figure {\ref{fig7}}.

\begin{figure}[h]
\begin{tikzpicture}[scale=.25,mine/.style={font=\itshape\large}]
\draw[thick]  [->] (-7,0) -- (-7,8);
\draw[thick]  [->] (-5,0) -- (-5,8);

\draw[thick]  [<->] (-4,4) -- (-1,4);
\draw (-2.5,5) node[mine] {$\Omega_{2a}$};

\draw[thick]  [<-] (0,8) .. controls (4,4) and (4,4) .. (0,0);
\draw[thick]  [<-] (4,8) .. controls (2.5,6.7) and (2,6)  .. (2.5,6.5);
\draw[thick]  (1.5,5.5) .. controls (0.5,4) and (0.5,4)  .. (1.5,2.5);
\draw[thick]   (4,0) .. controls (3,0.9)   .. (2.3,1.6);

\draw[thick]  [<-] (10,6) -- (18,6);
\draw[thick]    (10,8) -- (11.3,6.5);
\draw[thick]   [->] (12,5.5) -- (16,1);

\draw[thick] [<-] (18,8) -- (16.8,6.6);
\draw[thick]  (16.1,5.7) -- (14.5,3.7);
\draw[thick]  (13.8,3) -- (12,1);

\draw[thick]   [<->] (20,5) -- (24,5);
\draw[thick]  (22.5,6) node[mine] {$\Omega_{3a}$};

\draw[thick] [<-] (26,3) -- (34,3);
\draw[thick]  (27,8) -- (31.8,3.3);
\draw[thick]  [->]  (32.3,2.6) -- (34,1);
\draw[thick] [<-] (33,8) -- (30.4,5.2);
\draw[thick]  (29.7,4.4) -- (28.7,3.2);
\draw[thick]  (28.4,2.7) -- (27,1);
\end{tikzpicture}
\caption{Reidemeister move two and three.}
\label{fig8}
\end{figure}
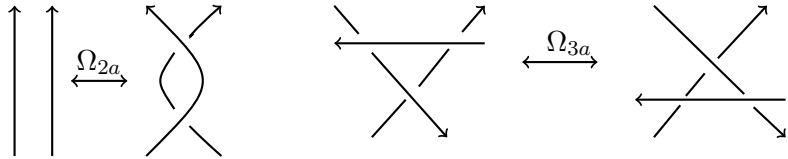

When we consider orientations, their are $8$ cases for the diagram inside $D_i$. When we consider coloring, each case has eight subcases. There will be too many cases. Fortunately, according to {\cite{P}}, we just need to consider one case $\Omega_{3a}$.  See figure {\ref{fig8}}.

To apply skein relations (\ref{eq1}) to $\Omega_{3a}$, we also need to color the link diagram. For each of the $D$ or $D'$, each of the three segments $a,b,c$ can be colored in 2 ways. So there are 8 cases. Let's first consider the case that they are all dotted. Then for $L$, crossings $1,2,3$ are of type $S_+,S_-,S_+$. Then for $L'$, crossings $1,2,3$ are of type $N_+,N_-,N_+$.

\begin{figure}[h]
\begin{tikzpicture}[scale=0.3,mine/.style={font=\itshape\large}]
\draw[thick]  [<-] (9,6) -- (19,6);
\draw[thick]   (10,8) -- (11.3,6.5);
\draw[thick]  [->] (12,5.5) -- (16,1);

\draw[thick][<-] (18,8) -- (16.8,6.6);
\draw[thick] (16.1,5.7) -- (14.5,3.7);
\draw[thick] (13.8,3) -- (12,1);
\draw[thick] (9.5,8.5) node {$a$};
\draw[thick] (9,7) node {$b$};
\draw[thick] (11.5,1) node {$c$};

\draw[thick] (15.5,3) node {$1$};p
\draw[thick] (12,7) node {$3$};
\draw[thick] (16,7) node {$2$};
\draw[thick] (14,-1) node {$D$};

\draw[thick][shift={(0,-.5)}]   [<->] (20,5) -- (24,5);
\draw[thick][shift={(0,-.5)}]  (22,6) node {$\Omega_{3a}$};

\draw[thick][shift={(0,1)}] [<-] (25,3) -- (35,3);
\draw[thick][shift={(0,1)}]  (27,8) -- (31.8,3.3);
\draw[thick][shift={(0,1)}]  [->]  (32.3,2.6) -- (34,1);
\draw[thick][shift={(0,1)}]  [<-] (33,8) -- (30.4,5.2);
\draw[thick][shift={(0,1)}]  (29.7,4.4) -- (28.7,3.2);
\draw[thick][shift={(0,1)}]  (28.4,2.7) -- (27,1);
\draw[thick][shift={(0,0)}]  (30,-1) node {$D'$};
\draw[thick] (26.5,8.5) node {$a$};
\draw[thick] (25,5) node {$b$};
\draw[thick] (26.3,2) node {$c$};

\draw[thick] (30,7) node {$1$};
\draw[thick] (29,3) node {$2$};
\draw[thick] (32,3) node {$3$};

\draw  (14,5) ellipse (4 and 4);
\draw[shift={(0,1)}]   (30,4) ellipse (4 and 4);
\end{tikzpicture}

\caption{Reidemeister move $\Omega_{3a}$.}
\label{fig9}
\end{figure}
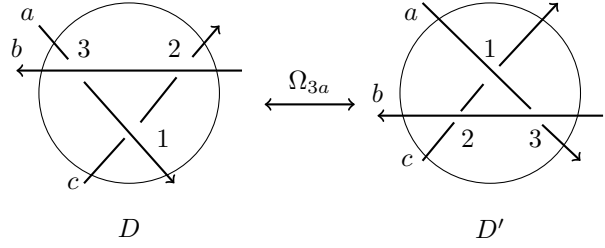

Now apply the skein relations to the diagrams. There are three crossings denoted by $1,2,3$ respectively. We get the following table (\ref{tab1}).

\begin{table}[ht]
\caption{Smooth inside crossings.}\label{tab1}
\begin{center}
\resizebox{12.5cm}{!}{\begin{tabular}[pos]{|l|l|l|l|l|l|l|l|l|l|}
\hline   & $a_1a_2a_3$ & $a_1a_2b_3$ & $a_1b_2a_3$ & $a_1b_2b_3$ & $b_1a_2a_3$ &$b_1a_2b_3$ & $b_1b_2a_3$ & $b_1b_2b_3$    \\
\hline $L$ & $B$ \tikz{\draw [scale=.1]  (0,0) ellipse (4 and 4);
\draw[thick] [scale=.1] (-1.5,3.6) arc (-215:35:1.9);
\draw [thick] [scale=.1] (-4,0) arc (120:-60:2);
\draw[thick] [scale=.1] (4,0) arc (60:240:2);
\draw [thick] [scale=.1] (0,-.5) ellipse (.8 and .8);}

& $A$ \tikz{\draw [scale=.1] (0,0) ellipse (4 and 4);
\draw [thick][scale=.1] (-1.5,3.6) arc (-215:35:1.9);
\draw [thick] [scale=.1](-4,0) arc (120:-60:2);
\draw [thick] [scale=.1](4,0) arc (60:240:2); } & $A$ \tikz{\draw [scale=.1] (0,0) ellipse (4 and 4);
\draw [thick][scale=.1] (-1.5,3.6) arc (-215:35:1.9);
\draw [thick] [scale=.1](-4,0) arc (120:-60:2);
\draw [thick] [scale=.1](4,0) arc (60:240:2); }

 & $E$  \tikz{\draw  [scale=.1] (0,0) ellipse (4 and 4);
\draw [thick][scale=.1] (-1.5,3.6) arc (-215:35:1.9);
\draw [thick][scale=.1] (-1.5,-3.6) arc (215:-35:1.9);
\draw [thick][scale=.1] (-4,0) -- (4,0);}

 & $A$ \tikz{\draw [scale=.1] (0,0) ellipse (4 and 4);
\draw [thick][scale=.1] (-1.5,3.6) arc (-215:35:1.9);
\draw [thick] [scale=.1](-4,0) arc (120:-60:2);
\draw [thick] [scale=.1](4,0) arc (60:240:2); }

& $F$  \tikz{\draw  [scale=.1] (0,0) ellipse (4 and 4);
\draw  [thick][scale=.1] (4,0) arc (300:120:2);
\draw [thick][scale=.1] (-4,0) arc (120:-60:2);
\draw  [thick] [scale=.1] (-1.5,3.6) -- (1.5,-3.6);}
& $G$ \tikz{\draw [scale=.1] (0,0) ellipse (4 and 4);
\draw [thick][scale=.1] (-4,0) arc (-120:60:2);
\draw [thick] [scale=.1](4,0) arc (60:240:2);
\draw [thick] [scale=.1] (1.5,3.6) -- (-1.5,-3.6);}

&

$C$  \tikz{\draw [scale=.1] (0,0) ellipse (4 and 4);
\draw [thick] [scale=.1](-4,0) arc (-120:60:2);
\draw [thick] [scale=.1](4,0) arc (300:120:2);
\draw  [thick][scale=.1] (-1.5,-3.6) arc (215:-35:1.9);}

\\
\hline $L'$ & $D$ \tikz{\draw [scale=.1]  (0,0) ellipse (4 and 4);
\draw [thick][scale=.1] (-4,0) arc (-120:60:2);
\draw [thick][scale=.1] (4,0) arc (300:120:2);
\draw [thick][scale=.1] (-1.5,-3.6) arc (215:-35:1.9);
\draw  [thick][scale=.1] (0,.5) ellipse (.8 and .8);}

& $C$ \tikz{\draw [scale=.1] (0,0) ellipse (4 and 4);
\draw [thick] [scale=.1](-4,0) arc (-120:60:2);
\draw [thick] [scale=.1](4,0) arc (300:120:2);
\draw  [thick][scale=.1] (-1.5,-3.6) arc (215:-35:1.9);}
& $C$  \tikz{\draw [scale=.1] (0,0) ellipse (4 and 4);
\draw [thick] [scale=.1](-4,0) arc (-120:60:2);
\draw [thick] [scale=.1](4,0) arc (300:120:2);
\draw  [thick][scale=.1] (-1.5,-3.6) arc (215:-35:1.9);}

& $E$  \tikz{\draw  [scale=.1] (0,0) ellipse (4 and 4);
\draw [thick][scale=.1] (-1.5,3.6) arc (-215:35:1.9);
\draw [thick][scale=.1] (-1.5,-3.6) arc (215:-35:1.9);
\draw [thick][scale=.1] (-4,0) -- (4,0);}

& $C$ \tikz{\draw [scale=.1] (0,0) ellipse (4 and 4);
\draw [thick] [scale=.1](-4,0) arc (-120:60:2);
\draw [thick] [scale=.1](4,0) arc (300:120:2);
\draw  [thick][scale=.1] (-1.5,-3.6) arc (215:-35:1.9);}

 & $F$  \tikz{\draw  [scale=.1] (0,0) ellipse (4 and 4);
\draw  [thick][scale=.1] (4,0) arc (300:120:2);
\draw [thick][scale=.1] (-4,0) arc (120:-60:2);
\draw  [thick] [scale=.1] (-1.5,3.6) -- (1.5,-3.6);}

& $G$ \tikz{\draw [scale=.1] (0,0) ellipse (4 and 4);
\draw [thick][scale=.1] (-4,0) arc (-120:60:2);
\draw [thick] [scale=.1](4,0) arc (60:240:2);
\draw [thick] [scale=.1] (1.5,3.6) -- (-1.5,-3.6);}

&  $A$ \tikz{\draw [scale=.1] (0,0) ellipse (4 and 4);
\draw [thick][scale=.1] (-1.5,3.6) arc (-215:35:1.9);
\draw [thick] [scale=.1](-4,0) arc (120:-60:2);
\draw [thick] [scale=.1](4,0) arc (60:240:2); } \\
\hline
\end{tabular}}
\end{center}
\end{table}

The table (\ref{tab1}) means the following. For example, the second row the third column is the result of $A$-type smoothing for crossings 1 of $L$ inside $D$. $a_2$ means that the second crossing uses $A$ type smoothing. $b_3$ will mean that the 3rt crossing uses $B$ type smoothing. This will be called using $a_1a_2b_3$-type smoothing.

Let $D_i$ be the diagram that inside a disk it is the same as $L$ in $D$, outside the disk it is the same as $Out (i)$. Let $D_i'$ be the diagram that inside a disk it is the same as $L'$ in $D'$, outside the disk it is the same as $Out (i)$.  If we smooth all outside crossings, we will have the following result. We have the following table (\ref{tab2}).

\begin{table}[ht]
\caption{Number of components of smoothings of $D_i$ and $D_i'$}\label{tab2}
\begin{center}
\resizebox{12.5cm}{!}{\begin{tabular}[pos]{|l|l|l|l|l|l|l|l|l|l|l|l|}
\hline $D_i' / D_i $ & $a_1a_2a_3$ & $a_1a_2b_3$ & $a_1b_2a_3$ & $a_1b_2b_3$ & $b_1a_2a_3$ &$b_1a_2b_3$ & $b_1b_2a_3$ & $b_1b_2b_3$    \\
\hline $Out (1)$ & 4/2 & 3/1 & 3/1  & 2/2 & 3/1   & 2/2  &  2/2 &  1/3 \\
\hline $Out (2)$ & 2/4 & 1/3 & 1/3  & 2/2 & 1/3   & 2/2  &  2/2 &  3/1 \\
\hline $Out (3)$ & 3/3 & 2/2 & 2/2  & 3/3 & 2/2   & 1/1  &  1/1 &  2/2 \\
\hline $Out (4)$ & 3/3 & 2/2 & 2/2  & 1/1 & 2/2   & 1/1  &  3/3 &  2/2 \\
\hline $Out (5)$ & 3/3 & 2/2 & 2/2  & 1/1 & 2/2   & 3/3   & 1/1  &  2/2 \\
\hline
\end{tabular}}
\end{center}
\end{table}

This is the table (\ref{tab2}) of number of components for smoothings of $D_i$ and $D_i'$. For example, the second row the second column is the result of $A$-type smoothing to all crossings of $D_1$. 4 means there are 4 disjoint circles. 2 means that if one applies $A$-type smoothing to all crossings of $D_1'$, there are 2 disjoint circles. For example, the third row the fourth column is $1/3$. This means that if one applies $a_1b_2a_3$-type smoothing to all crossings of $D_1$, there are 4 disjoint circles. 2 means that if one apply $a_1b_2a_3$-type smoothing to all crossings of $D_1'$, there are 2 disjoint circles.

To get $<D_i>=<D_i'>$, the second row of table (\ref{tab2}) gives the following equation.

$a_na_n'a_nd^4+  a_na_n'b_nd^3+   a_nb_n'a_nd^3+  a_nb_n'b_nd^2+   b_na_n'a_nd^3+  b_na_n'b_nd^2+  b_nb_n'a_nd^2+  b_nb_n'b_nd= \\
a_sa_s'a_sd^2+  a_sa_s'b_sd+     a_sb_s'a_sd+    a_sb_s'b_sd^2+   b_sa_s'a_sd+    b_sa_s'b_sd^2+  b_sb_s'a_sd^2+  b_sb_s'b_sd^3.$

The above is the equation from $Out (1)$. Denote $x_1=a_na_n'a_n, x_2=a_na_n'b_n, \cdots, x_8=b_nb_n'b_n, y_1=a_sa_s'a_s, \cdots ,y_8=b_sb_s'b_s$, we get a linear equation for variables $x_1,\cdots ,y_8$.
The following is the coefficient matrix/table of the system of equations from $Out (1)$ -- $Out (5)$.

\begin{table}[ht]
\caption{Coefficient matrix of the system of equations}\label{tab3}
\begin{center}
\resizebox{12.5cm}{!}{\begin{tabular}[pos]{|l|l|l|l|l|l|l|l|l|l|l|l|l|l|l|l|l|l|l|l|l|l|l|}
\hline $       $ & $x_1$ & $x_2$ & $x_3$  & $x_4$ & $x_4$   & $x_6$  &  $x_7$ &  $x_8$ & $y_1$ & $y_2$ & $y_3$  & $y_4$ & $y_5$  & $y_6$  &  $y_7$ &  $y_8$ \\
\hline $Out (1)$ & $d^4$ & $d^3$ & $d^3$  & $d^2$ & $d^3$   & $d^2$  &  $d^2$ &  $d$ & $d^2$ & $d$ & $d$  & $d^2$ & $d$  & $d^2$  &  $d^2$ &  $d^3$ \\
\hline $Out (2)$ & $d^2$ & $d$ & $d$  & $d^2$ & $d$   & $d^2$  &  $d^2$ & $ d^3$ & $d^4$ & $d^3$ & $d^3$  & $d^2$ & $d^3$   & $d^2$  &  $d^2$ &  $d $\\
\hline $Out (3)$ & $d^3$ & $d^2$ & $d^2 $ & $d^3$ & $d^2 $  & $d$  &  $d$ &  $d^2$ & $d^3$ & $d^2$ & $d^2$  & $d^3$ & $d^2$   & $d$  &  $d$ &  $d^2$  \\
\hline $Out (4)$ & $d^3$ & $d^2 $& $d^2$  & $d $  & $d^2$   & $d$  &  $d^3$ &  $d^2$ & $d^3$ & $d^2$ & $d^2$  & $d $  & $d^2$   & $d $ &  $d^3$ &  $d^2$ \\
\hline $Out (5)$ & $d^3$ & $d^2$ & $d^2 $ & $d $ & $d^2$  & $d^3 $  & $d$  & $ d^2$ &  $d^3$ & $d^2$ & $d^2 $ & $d$  & $d^2$  &$ d^3 $  & $d $ & $ d^2$ \\
\hline
\end{tabular}}
\end{center}
\end{table}

Notice that the coefficients of $y_1,\cdots ,y_8$ should have a $-$ sign. We do not put the negative sign there, one can regard this as that $y_1,\cdots ,y_8$ and their coefficients are on the righthand side of the equation.

Let $Eq_i$ denote the equation corresponding to $Out (i)$. Then from $Eq_5-Eq_4$ and $Eq_4$ we get $x_6-x_7 =y_6-y_7=\alpha$. From $Eq_3-Eq_4$ and $Eq_4$  we get $x_4-x_7 =y_4-y_7=\beta$.
From $Eq_1-Eq_4$ and $Eq_4$ we get $dx_7+x_8 =dy_1+y_2+y_3+y_5+dy_7=\gamma$.
From $Eq_2-Eq_4$ and $Eq_4$ we get $dy_7+y_8 =dx_1+x_2+x_3+x_5+dx_7=\delta$.
Plug those into $Eq_4$ we get $(2-d^2)(x_7-y_7)=0$.
If we take the solution $x_7=y_7$, then $x_6=y_6$ and $x_4=y_4$, $x_8=dy_1+y_2+y_3+y_5$, $y_8=dx_1+x_2+x_3+x_5$.

In other words, we have the following.

\begin{equation}\label{eq3}
 \left\{
\begin{lgathered}
b_nb_n'a_n=b_sb_s'a_s,  \ \ \ \ \  b_na_n'b_n=b_sa_s'b_s,  \ \ \ \ \  a_nb_n'b_n=a_sb_s'b_s  \\
b_nb_n'b_n =da_sa_s'a_s+a_sa_s'b_s+a_sb_s'a_s+b_sa_s'a_s\\
b_sb_s'b_s =da_na_n'a_n+a_na_n'b_n+a_nb_n'a_n+b_na_n'a_n\\
\end{lgathered}
\right.
\end{equation}

For each of the $D$ or $D'$, if we change the color of the three segments $a,b,c$, we shall get the following results. In table  (\ref{tab3}) line 1, $a1b2c1$ means that segment $a$ choose color 1 (solid), segment $b$ choose color 2 (dotted),   segment $c$ choose color 1 (solid). In that column, $ N_+W_-E_+ $ means that the three crossings $1,2,3$ of $L$ are of type $N_+,W_-,E_+$ respectively, and the three crossings $1,2,3$ of $L'$ are of type $S_+,E_-,W_+$ respectively.

\begin{table}[ht]
\caption{Number of components of smoothings of $D_i$ and $D_i'$}\label{tab4}
\begin{center}
\resizebox{12.5cm}{!}{\begin{tabular}[pos]{|l|l|l|l|l|l|l|l|l|l|l|l|}
\hline  & $a1b1c1$ & $a1b1c2$ & $a1b2c1$ & $a1b2c2$ & $a2b1c1$ &$a2b1c2$ & $a2b2c1$ & $a2b2c2$    \\
\hline $L$ & $N_+N_-N_+ $ & $ W_+E_-N_+ $ & $ N_+W_-E_+ $ & $ W_+S_-E_+ $ & $ E_+N_-W_+   $ & $ S_+E_-W_+  $ & $  E_+W_-S_+ $ & $ S_+S_-S_+$ \\
\hline $L'$ & $S_+S_-S_+ $ & $ E_+W_-S_+ $ & $S_+E_-W_+  $ & $ E_+N_-W_+ $ & $ W_+S_-E_+  $ & $ N_+W_-E_+  $ & $  W_+E_-N_+ $ & $  N_+N_-N_+$ \\
\hline
\end{tabular}}
\end{center}
\end{table}

If we change the color of $a,b,c$ to $a1b2c1$, then the above discuss is almost the same except we need to change the subscripts in (3) to $n,w',e$ and $s,e',w$.

\begin{equation}\label{eq4}
 \left\{
\begin{lgathered}
b_nb_w'a_e=b_sb_e'a_w,  \ \ \ \ \  b_na_w'b_e=b_sa_e'b_w,  \ \ \ \ \  a_nb_w'b_e=a_sb_e'b_w  \\
b_nb_w'b_e =da_sa_e'a_w+a_sa_e'b_w+a_sb_e'a_w+b_sa_e'a_w\\
b_sb_e'b_w =da_na_w'a_e+a_na_w'b_e+a_nb_w'a_e+b_na_w'a_e\\
\end{lgathered}
\right.
\end{equation}

Now, we go over all cases in Table  (\ref{tab4}), we will get the following equations from $\Omega_{3a}$.

\begin{equation}\label{R3}
 \left\{
\begin{lgathered}
b_nb_n'a_n=b_sb_s'a_s,  \ \ \ \ \  b_na_n'b_n=b_sa_s'b_s,  \ \ \ \ \  a_nb_n'b_n=a_sb_s'b_s  \\
b_nb_n'b_n =da_sa_s'a_s+a_sa_s'b_s+a_sb_s'a_s+b_sa_s'a_s\\
b_sb_s'b_s =da_na_n'a_n+a_na_n'b_n+a_nb_n'a_n+b_na_n'a_n\\
b_nb_w'a_e=b_sb_e'a_w,  \ \ \ \ \  b_na_w'b_e=b_sa_e'b_w,  \ \ \ \ \  a_nb_w'b_e=a_sb_e'b_w  \\
b_nb_w'b_e =da_sa_e'a_w+a_sa_e'b_w+a_sb_e'a_w+b_sa_e'a_w\\
b_sb_e'b_w =da_na_w'a_e+a_na_w'b_e+a_nb_w'a_e+b_na_w'a_e\\
b_wb_s'a_e=b_eb_n'a_w,  \ \ \ \ \  b_wa_s'b_e=b_ea_n'b_w,  \ \ \ \ \  a_wb_s'b_e=a_eb_n'b_w  \\
b_wb_s'b_e =da_ea_n'a_w+a_ea_n'b_w+a_eb_n'a_w+b_ea_n'a_w\\
b_eb_n'b_w =da_wa_s'a_e+a_wa_s'b_e+a_wb_s'a_e+b_wa_s'a_e\\
b_wb_e'a_n=b_eb_w'a_s,  \ \ \ \ \  b_wa_e'b_n=b_ea_w'b_s,  \ \ \ \ \  a_wb_e'b_n=a_eb_w'b_s  \\
b_wb_e'b_n =da_ea_w'a_s+a_ea_w'b_s+a_eb_w'a_s+b_ea_w'a_s\\
b_eb_w'b_s =da_wa_e'a_n+a_wa_e'b_n+a_wb_e'a_n+b_wa_e'a_n\\
\end{lgathered}
\right.
\end{equation}

\bigskip
\subsection{Simplications}

The equations (\ref{R2}) from Reidemeister II can be simplified to get the following equations.

\begin{equation}\label{RS2}
 \left\{
\begin{lgathered}
a_e'=a_w^{-1}, a_w'=a_e^{-1}, a_s'=a_s^{-1},a_n'= a_n^{-1}\\
b_e'=b_w^{-1}, b_w'=b_e^{-1}, b_s'=b_s^{-1},b_n'= b_n^{-1}\\
-d=\frac{a_w}{b_w}+\frac{b_w}{a_w}=\frac{a_e}{b_e}+\frac{b_e}{a_e}=\frac{a_s}{b_s}+\frac{b_s}{a_s}=\frac{a_n}{b_n}+\frac{b_n}{a_n}\\
\end{lgathered}
\right.
\end{equation}

One solution is to introduce new variables $e,w,n,s,a,b$ and let $$a_n=na,b_n=nb,a_s=sa,b_s=sb,a_w=wa,b_w=wb,a_e=ea,b_e=eb, -d=\frac{a}{b}+\frac{b}{a}.$$ Then $$-d=\frac{a_w}{b_w}+\frac{b_w}{a_w}=\frac{a_e}{b_e}+\frac{b_e}{a_e}=\frac{a_s}{b_s}+\frac{b_s}{a_s}=\frac{a_n}{b_n}+\frac{b_n}{a_n}$$ is satisfied. Then $$a_e'=\frac{1}{wa}, a_w'=\frac{1}{ea}, a_s'=\frac{1}{sa},a_n'= \frac{1}{na}, b_e'=\frac{1}{wb}, b_w'=\frac{1}{eb}, b_s'=\frac{1}{sb},b_n'= \frac{1}{nb}.$$

Plug those into equation \label{R3}, we get one extra equation: $n=s$. Then all equations \label{R3} are satisfied.

In summary, we have free variables $w,e,n,a,b$. The other variables can be derived from them.

\begin{equation}\label{All}
 \left\{
\begin{lgathered}
a_n=a_s=na,b_n=b_s=nb,a_w=wa,b_w=wb,a_e=ea,b_e=eb,\\
a_n'=a_s'= \frac{1}{na},b_n'=b_s'= \frac{1}{nb},a_w'=\frac{1}{ea}, b_w'=\frac{1}{eb}, a_e'=\frac{1}{wa},b_e'=\frac{1}{wb}\\
d=-\frac{a}{b}-\frac{b}{a}\\
\end{lgathered}
\right.
\end{equation}

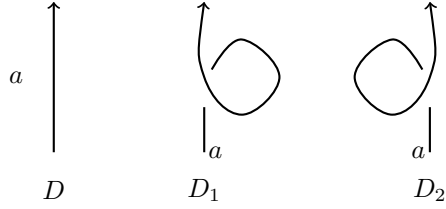
\begin{figure}[h]
\begin{tikzpicture} [scale=.5]
\draw [thick] [<-] (-1,2) -- (-1,-2);
\draw [thick] (-2,0) node {$a$};
\draw [thick] (-1,-3) node {$D$};

\draw [thick] [<-] plot[smooth, tension=.7] coordinates {(3,2) (3,0) (4,-1) (5,0) (4,1) (3.2,0.2)};
\draw [thick] plot[smooth, tension=.7] coordinates {(3,-0.8) (3,-2)};
\draw [thick] (3.3,-2) node {$a$};
\draw [thick] (3,-3) node {$D_1$};

\draw [thick] [<-] plot[smooth, tension=.7] coordinates {(9,2) (9,0) (8,-1) (7,0) (8,1) (8.8,0.2)};
\draw [thick]  plot[smooth, tension=.7] coordinates {(9,-0.8) (9,-2)};
\draw [thick] (8.7,-2) node {$a$};
\draw [thick] (9,-3) node {$D_2$};
\end{tikzpicture}
\caption{Reidemeister move one invariance.}
\label{fig10}
\end{figure}

If we consider Reidemeister move I (Figure {\ref{fig10}}), choose one color for the arc $a$, say, dotted, then the crossing type of $D_1$ is $N_-$. we have $$<D_1>=(da_n'+b_n')<D>=(-(\frac{a}{b}+\frac{b}{a})\frac{1}{na}+\frac{1}{nb})<D>=-\frac{b}{na^2}<D>.$$

The crossing type of $D_2$ is $S_+$. we have $$<D_1>=(da_n+b_n)<D>=(-(\frac{a}{b}+\frac{b}{a})na+nb)<D>=-\frac{na^2}{b}<D>.$$

If we choose one color for the arc $a$, say, solid, then the crossing type of $D_1$ is $S_-$. we have $$<D_1>=(da_n'+b_n')<D>=(-(\frac{a}{b}+\frac{b}{a})\frac{1}{na}+\frac{1}{nb})<D>=-\frac{b}{na^2}<D>.$$

The crossing type of $D_2$ is $N_+$. we have $$<D_1>=(da_n+b_n)<D>=(-(\frac{a}{b}+\frac{b}{a})na+nb)<D>=-\frac{na^2}{b}<D>.$$

Let $W(D)$ denote the writhe of the diagram $D$. Like the construction of Kauffman bracket, let $F(D)=(-\frac{b}{na^2})^{W(D)}<D>$. Then $F(D)$ is invariant under Reidemeister moves. However, it depend on coloring of the link diagram. For a knot diagram, there are only two different colorings. If we change the coloring, we shall get $\overline{F}(D)$. It can be obtained from $F(D)$ as follows. Define $\overline{e}=w,\overline{w}=e$. Then we have

\begin{equation}\label{All}
 \left\{
\begin{lgathered}
\overline{a_n}=\overline{a_s}=na,\overline{b_n}=\overline{b_s}=nb,\overline{a_w}=ea,\overline{b_w}=eb,\overline{a_e}=wa,\overline{b_e}=wb,\\
\overline{a_n}'=\overline{a_s}'= \frac{1}{na},\overline{b_n}'=\overline{b_s}'= \frac{1}{nb},\overline{a_w}'=\frac{1}{wa}, \overline{b_w}'=\frac{1}{wb}, \overline{a_e}'=\frac{1}{ea},\overline{b_e}'=\frac{1}{eb}\\
\overline{d}=-\frac{a}{b}-\frac{b}{a}\\
\end{lgathered}
\right.
\end{equation}.
If we change each variable $x$ to $\overline{x}$ , $F(D)$ becomes  $\overline{F}(D)$.

In general, For a link $L$, choose one link diagram $D$, let $\Lambda$ be the set of all colorings. If one choose one coloring $\lambda \in \Lambda$, one will get $F(D,\lambda )$. Now we have the following theorem.
\begin{theorem}\label{sec:6}
Using the skein relations (1), if the variables satisfies equation (7), then $\{F(D),\overline{F}(D)\}$ is a knot invariant.

For a link $L$, choose one link diagram $D$, then $F(L) =\{F(D,\lambda ) \mid  \lambda \in \Lambda\}$ is a multiple-valued link invariant.
\upshape
\end{theorem}

The invariant $F(L)$ is also an virtual knot invariant. Given an virtual link diagram, we can also color it as follows. Choose two colors, say, solid and dotted. The crossing points divide any link component into even number of segments. Pick one segment and assign one color to it. Then change to another color whenever pass one crossing point. But the color does not change when pass one virtual crossing point. Do this to each component, we  get a bicolor link diagram.
For crossings in a virtual link diagram, we also use the skein relations (1), but we do not apply skein relations to virtual crossings. For a diagram $D$ contains only virtual crossings, if it has $n$ link component, let $<D>=d^n$. The the above theorem is valid.

\bigskip\noindent {\bf Example: A one variable inavriant.}

S. Nelson, M. Orrison, V. Rivera {\cite{S}} modified Kauffman's construction of bracket. Their invariant $\Phi^{\beta}_X$ takes value in a finite ring $Z_2[t]/(1+t+t^3)$. Their choice of the coefficients are as follows.

\begin{equation}\label{N}
 \left\{
\begin{lgathered}
a_n=a_s=1,b_n=b_s=t,a_w=1+t^2,b_w=1,a_e=1+t,b_e=t+t^2,\\
a_n'=a_s'= 1,b_n'=b_s'= 1+t^2,a_w'=t+t^2, b_w'=1+t, a_e'=t,b_e'=1\\
d=1+t+t^2\\
\end{lgathered}
\right.
\end{equation}

The coefficients lie in $Z_2[t]/(1+t+t^3)$. We can lift them to $Z[t,t^{-1}]$ as follows.

\begin{equation}\label{Y}
 \left\{
\begin{lgathered}
a_n=a_s=1,b_n=b_s=t,a_w=1+t^2,b_w=t(1+t^2),a_e=1+t,b_e=t(1+t),\\
a_n'=a_s'= 1,b_n'=b_s'= 1/t,a_w'=\frac{1}{1+t}, b_w'=\frac{1}{t(1+t)}, a_e'=\frac{1}{1+t^2},b_e'=\frac{1}{t(1+t^2)}\\
d=-t-\frac{1}{t}\\
\end{lgathered}
\right.
\end{equation}

For this choice, $a=1,b=t,w=1+t^2,e=1+t,n=s=1$.

\begin{flalign*}
N_+ = H+ t V, &\ \ \ \ \   N_- = H+ (1+t^2) V,  \\
S_+ = H+ t V, &\ \ \ \ \     S_- = H+ (1+t^2) V,  \\
E_+ =(1+t) H+ (t+t^2) V,  &\ \ \ \ \   E_- =t H+ V,   \\
W_+ =(1+t^2) H+  V,  &\ \ \ \ \    W_- =(t+t^2)  H+ (1+t)  V,   \\
<\cup^n \tikz{\draw (0,0) circle (0.12cm);}> &=(1+t+t^2)^n.
\end{flalign*}

This invariant, however, does not gives new results for classical knots. One can easily find out that for any knot diagram, using bicolor, one can only get $N_+,N_-,S_+,S_-$ type crossings. If one extend to virtual knots, this invariant turns out to give more information than the Jones polynomial.  

\section{Tricoloring invariant}
If we use three colors, we will get a nontrivial results for knots. 

Pick up three colors, say, red, blue and green. A link diagram can be colored with the following rules:
At each crossing, either all three colors are present or only one color is present. If one uses only one color we say that it is a trivial
tricoloring. The number of different tricolorings (trivial cloring is allowed) is denoted by $tri(D)$.

\begin{lem}{\cite{P}} $tri(L)$ is always a power of 3.
\end{lem}

\begin{theorem}{\cite{P}}
(a) $tri(L) = 3|V_L^2(e^{2\pi i/6})|$\\
(b) $tri(L) = 3|FL(1, -1)|$\\
\end{theorem}

Here $V_L$ is the Jones polynomial, $FL$ is the two-variable Kauffman polynomial. For example, Figure eight knot only has trivial tricoloring, then $tri(L) = 3$. $V(7_4)=t-2*t^2+ 3*t^3-2*t^4+ 3*t^5-2*t^6+ t^7-t^8$, then $tri (7_4)=9$. This means that $7_4$ has only one nontrivial coloring up to permutation of the colors.

For tricolored link diagram, we can define the following skein relations. See Fig. ~\ref{fig1}. If the the three arcs have same color 

$$<L_+>=x<H>+x^{-1}<V>,<L_->=x^{-1}<H>+x<V>.$$

Otherwise, $$<L_+>=y<H>+y^{-1}<V>, <L_->=y^{-1}<H>+y<V>.$$

$$<D\  \tikz{\draw [thick][scale=.25] (0,1) -- (0,0);
\draw [thick][scale=.25] (0,0) -- (1,0);
\draw [thick][scale=.25] (1,0) -- (1,1);}\  \tikz{\draw (0,0) circle (0.14cm);}> =(-x^2-x^{-2}) <D>$$

One can easily check that this bracket is invariant under Reidemeister moves II and III. For Reidemeister moves I, the three arcs always have same color. Hence if we let $V(D)=(-x^3)^{-w(D)}<D>$, we shall get a knot invariant.

\begin{theorem}
$V(D)=(-x^3)^{-w(D)}<D>$ is a two variable knot invariant. 
\end{theorem}

\bigskip\noindent {\bf Example: the $7_4$ knot.}

As showed before, $tri (7_4)=9$. Hence $7_4$ has only one nontrivial coloring up to permutation of the colors. If one take any $7$ crossing diagram of $7_4$ and color it, then one will find that there is one crossing with three arcs having same color, the other six crossings with three arcs having different color.  

On the other hand, $7_4$ is alternating, hence the bracket is ``faithful". Therefor, for any diagram, any nontrivial tricoloring, there exists one crossing with same color, and at least 6 crossings with different colors.

\section*{Acknowledgements}
The author would like to thank Ruifeng Qiu, Jiajun Wang, Ying Zhang, Xuezhi Zhao for helpful discussions.

\end{document}